\documentclass[11pt,a4paper,french]{amsart}

\usepackage{babel}
\language1 \theoremstyle{plain}
\newtheorem{theorem}{Th\'eor\`eme}[section]
\newtheorem{lemma}[theorem]{Lemme}
\newtheorem{proposition}[theorem]{Proposition}
\newtheorem{corollary}[theorem]{Corollaire}

\theoremstyle{definition}

\theoremstyle{remark}
\newtheorem*{remark}{Remarque}

\allowdisplaybreaks

\DeclareMathOperator{\real}{Re}

\title[Exemples de points extr\'emaux]{Le c\^one des fonctions
plurisousharmoniques n\'egatives et une conjecture de Coman}
\dedicatory{D\'edi\'e a J\'ozef Siciak}
\author{Magnus Carlehed et Jan Wiegerinck}
\email{magnus.carlehed@foreningssparbanken.se,janwieg@wins.uva.nl}
\address{Magnus Carlehed, Renstiernas gata 14, 116 28 Stockholm, Sweden,
\newline Jan Wiegerinck, Korteweg-de Vries Institute for
 Mathematics, University of Amsterdam,
Plantage Muidergracht 24, 1018 TV, Amsterdam, The Netherlands}

\begin{document}

\begin{abstract}
Les fonctions plurisousharmoniques n\'egatives dans un domaine
$\Omega$ de $\mathbb{C}^n$ forment un c\^one convexe. Nous
consid\'erons les points extr\'e\-maux de ce c\^one, et donnons
trois exemples. En particulier, nous traitons le cas de la
fonction de Green pluricomplexe. Nous calculons celle du bidisque,
lorsque les p\^oles se situent sur un axe. Nous montrons que cette
fonction ne se confonde pas avec la fonction de Lempert
correspondante. Cela donne un contre-exemple \`a une conjecture de
Dan Coman.
\end{abstract}

\maketitle

\section{Introduction}

Cet article a deux intentions : de donner quelques exemples de
points extr\'emaux dans le\,c\^one des fonctions
plurisousharmoniques n\'egatives, et de donner des exemples de
domaines convexes o\`u la fonction de Green pluricomplexe \`a
plusieurs p\^oles ne se confonde pas avec la fonction de Lempert
correspondante. Les deux buts sont li\'es. Les r\'esultats de l'
article actuel sont signal\'es dans les pr\'epublications
\cite{us1} et \cite{us2}.

Soit $C$ un c\^one convexe de sommet 0 dans un espace vectoriel
$V$. Un point $x\in C$ est contenu dans une g\'en\'eratrice
extr\'emale si $x=x_1+x_2$ o\`u $x_1,x_2\in C$ entra\^\i ne
$x_1=\lambda_1 x$, $x_2=\lambda_2 x$ , o\`u $\lambda_1,\lambda_2\ge
0$. Par abus de langage on appelle $x$ un point extr\'emal, voir
\cite{cho}. Si l'on munit $C$ d'une topologie m\'etrisable, le
sous-ensemble $E$ des points extr\'emaux sera un ensemble
$G_\delta$. On appelle $C$ {\em c\^one \`a base compact}, s'il
existe un hyperplan ferm\'e $H$ et un compact $K\subset H$ telle
que $C=\{tx\ : \, x\in K, t>0\}$. Maintenant, on a le th\'eor\`eme
suivant, \cite{cho} p. 140.
\begin{theorem}(Choquet) Soit $C$ un c\^one convexe \`a base compact, m\'etri\-sable,
alors pour chaque $q\in C$ il existe une mesure de probabilit\'e
$\mu_q$, concentr\'ee sur $E$, telle que $$f(q)=\int_C f(x)
d\mu_q(x)$$ pour toute fonction $f\in V'$.
\end{theorem}
En cons\'equence, il est int\'eressant de caract\'eriser $E$.

Dans cet article nous consid\'erons le c\^one des fonctions
plurisousharmoniques n\'egatives dans un domaine $\Omega$ de
$\mathbb{C}^n$, avec la topologie induite par celle de
$L_{loc}^1(\Omega)$. C'est connu qu'il est un c\^one \`a base
compact: pour $H$ on prend $\{ f\in L_{loc}^1(\Omega):\, \int_G
f\, dV=-1\}$, o\`u $G$ est un ouvert quelconque, relativement
compacte, de $\Omega$. Voir \cite{hor} p. 149 et 229. Nous dirons
qu'une fonction $u$ plurisousharmonique n\'egative est {\it
extr\'emale} si elle est extr\'emale au sens de c\^one convexe et
{\it maximale} si $(dd^c u)^n=0$. Pour un compact $K$ dans
$\Omega$ on d\'efinit {\it la fonction extr\'emale relative} comme
l'enveloppe sup\'erieure de la famille $U=\{v\in PSH(\Omega):v\leq
0,v|K\leq -1\}$. Remarquons qu'a priori cette fonction peut \^etre
non-extr\'emale au sens convexe, malgr\'e le langage.

La {\em fonction de Green pluricomplexe} (\`a plusieurs p\^oles)
appartient au c\^one des fonctions plurisousharmoniques
n\'egatives. Elle est introduite par Lelong \cite{Lelong} ainsi :
Soit $\Omega$ un domaine de $\mathbb{C}^n$, et
$$A=\{(w_1,\nu_1),\ldots,(w_k,\nu_k)\}\subset\Omega\times\mathbb{R}^+,$$
o\`u $\mathbb{R}^+=[0,\infty)$. On appelle $w_1,\ldots,w_k$ les
{\em p\^oles} et $\nu_1,\ldots,\nu_k$ les {\em poids}. On pose
$$U_{A,\Omega}=\{u\in
PSH(\Omega):u(\zeta)-\nu_j\log|\zeta-w_j|\leq C_u, \zeta\to
w_j,j=1,\ldots,k\},$$ et on d\'efinie la fonction de Green
pluricomplexe $$g(z,A)=g_\Omega(z,A)=\sup\{u(z):u\in
U_{A,\Omega},u\leq 0\}.$$ S'il s'agit d'un seul p\^ole $w$ et le
poids est \'egal \`a 1, on \'ecrit normalement $g(z,w)$.

Le probl\`eme d'extr\'emalit\'e a \'et\'e \'etudi\'e par Cegrell
et Thor\-bi\"ornson \cite{johan}. Soient d'abord $n=1$ et
$\Omega=D$, le disque unit\'e. Dans ce cas ils ont montr\'e qu'une
fonction sousharmonique n\'egative $\varphi$ est extr\'emale si et
seulement si, ou bien
$$\varphi(z)=k\log|\frac{z-z_0}{1-z\bar{z_0}}|,\quad k\geq 0,\quad
z_0\in D,$$ ou bien $$\varphi(z)=kP(z,\xi_0),\quad k\leq
0,\quad\xi_0\in\partial D,$$ o\`u on note $P$ le noyau de Poisson
de $D$.

Passons au cas $n\geq 2$. Les m\^emes auteurs ont g\'en\'eralis\'e
leur pr\'ec\'edent r\'esultat, en montrant que dans un domaine
$\Omega\subset\mathbb{C}^n$, la fonction de Green pluricomplexe
{\em \`a un seul p\^ole} $g(\cdot,w)$ est toujours extr\'emale
pour chaque $w\in\Omega$. Si celles-ci donnaient toutes les
fonctions extr\'emales qui s'annulent au bord, alors le
th\'eor\`eme de Choquet 
donnerait que toute fonction plurisousharmonique $u$ n\'egative
ayant $0$ pour valeur au bord
 pourrait
\^etre repr\'esent\'ee de la fa\c con suivante :
$$u(z)=\int_\Omega g(z,w) d\mu_u(w),$$ o\`u $\mu_u$ est une mesure
positive. Ce type de potentiels a \'et\'e \'etudi\'e par le
premier auteur \cite{magnus}, qui a montr\'e que la masse de
Monge-Amp\`ere d'un potentiel born\'e dans la boule est toujours
une mesure absolument continue. Ceci sugg\`ere que les fonctions
de Green forment un tr\`es petit sous-ensemble de $E$.

Ici, nous montrons que si la fonction de Green \`a plusieurs
p\^oles est d\'ecompos\'ee, alors les deux parties sont
\'egalement des fonctions de Green ayant les m\^emes p\^oles, mais
peut-\^etre d'autres poids. En particulier, les fonctions
suivantes sont extr\'emales : 1) la fonction de Green avec deux
p\^oles de poids 1 dans la boule unit\'e, 2) la fonction de Green
avec deux p\^oles $(a,0)$ et $(b,0)$ de poids 1 dans le bidisque
de $\mathbb{C}^2$. Nous profitons du fait que ces fonctions sont
connues explicitement. Au vu de cela on pourrait se demander si
des fonctions n\'egatives, et maximales dans un ensemble assez
grand, sont extr\'emales. Il est un peu surprenant que dans le cas
du bidisque, les fonctions de Green avec des poids distincts, ne
soient pas extr\'emales. Nous montrons cela par un calcul
explicite de ces
 fonctions.

Ensuite, nous montrons que la fonction $\max\{\log|z|,-1\}$ dans la boule
unit\'e (la fonction
extr\'emale relative d'une boule plus petite) est extr\'emale. Il est
tentant de conjecturer que
toute fonction extr\'emale relative d'un compact assez r\'egulier
appartient aussi \`a $E$.

Toutes les exemples donn\'es jusqu'ici concernent des fonctions
qui s'annulent au bord du domaine. Nous discutons bri\`e\-vement
des fonctions qui ne s'annulent pas sur tout le bord, et nous
donnons un exemple \'el\'ementaire.

Passons \`a notre deuxi\`eme sujet. On introduit d'abord la
fonction de Lempert. Soit $D$ le disque unit\'e, et $\Omega$ et
$A$ comme plus haut. Pour tout point $z\in\Omega$, on note
$F_z=F_{z,A}$ la famille des application analytiques
$f:D\to\Omega$, telles que $f(0)=z$ et qu'il existe des points
$\zeta_1,\ldots,\zeta_k\in D$ avec $f(\zeta_j)=w_j$,
$j=1,\ldots,k$. On appelle ces applications, aussi bien que leurs
images, des {\em disques analytiques ajust\'es \`a $A$}. Pour tout
$f\in F_z$ on d\'efinie $d(f)=\sum_{j=1}^k \nu_j\log|\zeta_j|$ et
la {\em fonction de Lempert}
$\delta(z,A)=\delta_\Omega(z,A)=\inf\{d(f):f\in F_z\}$. On voit
facilement que $\delta(\cdot,A)\geq g(\cdot,A)$ avec \'egalit\'e
si et seulement si $\delta$ est plurisousharmonique. Un
th\'eor\`eme remarquable de Lempert \cite{Lempert} dit que cela
est vraiment le cas si $\Omega$ est convexe et $k=1$. Beaucoup
plus tard Coman \cite{dan} a montr\'e que c'est aussi le cas si
$\Omega$ est la boule unit\'e, $k=2$, et les poids sont \'egaux.

On d\'efinit $$\delta^A(z)=\min_{\emptyset\neq B\subseteq
A}\delta(z,B).$$ \'Evidemment, on a $$\delta(\cdot,A)\geq
\delta^A(\cdot)\geq g(\cdot,A).$$ Coman \cite{dan} a conjectur\'e
que, dans des domaines convexes born\'es, la deuxi\`eme
in\'egalit\'e est toujours une \'egalit\'e. Plus tard, Wikstr\"om
\cite{Frank}, 
a montr\'e que dans ces domaines, la
premi\`ere in\'egalit\'e est toujours une \'egalit\'e. En
cons\'equence, il reformule la conjecture sous le forme
$\delta(\cdot,A)=g(\cdot,A).$

Wikstr\"om a aussi montr\'e le th\'eor\`eme suivant.

\begin{theorem}\cite{Frank}{\em \ corollaire 2.3.}\label{Frank1}
Soit $\Omega$ un domaine born\'e et taut de $\mathbb{C}^n$, et soit $A$
comme plus haut. Alors pour
tout $z\in\Omega$ il existe un disque analytique $f$, tel que $f(0)=z$,
passant par un
sous-ensemble (non-vide) $\{w_{j_1},\ldots,w_{j_m}\}$ de
$\{w_1,\ldots,w_k\}$, tel que $d(f)$ est
\'egal \`a la borne inf\'erieure dans la d\'efinition de
$$\delta^A(z)=\delta(z,\{(w_{j_1},\nu_{j_1}),\ldots,(w_{j_m},\nu_{j_m})\}).$$
\end{theorem}

Motiv\'e par ce th\'eor\`eme, il propose la d\'efinition suivante.
Si $\delta(z,A)=d(f)$ pour quelque disque $f$ passant par $z$ et
par un sous-ensemble non-vide de $A$, alors $f$ est appel\'e un
{\em disque extr\'emal} pour $z$ et $A$. Ensuite il caract\'erise
les disques extr\'emaux dans des domains convexes, ce qui
g\'en\'eralise un th\'eor\`eme de Lempert. Il fait remarquer que
les poids sont invisibles dans la caract\'erisation, ce qui
indique que la conjecture peut \^etre probl\'ematique, ou en tous
cas impossible de montrer en utilisant cette m\'ethode.

Dans cet article nous donnons un contre-exemple \`a la conjecture de Coman.
En effet on a le
th\'eor\`eme suivant.

\begin{theorem}\label{counterexample}
Soit $\Omega$ le bidisque d'unit\'e et soit $A=\{((a,0),1),((b,0),2)\}$
avec $0<|a|,|b|<1$. Alors $\delta(z,A)\not\equiv g(z,A)$.
\end{theorem}

Il s'ensuivra un m\^eme r\'esultat pour des domaines convexes,
suffisamment proche de $\Omega$. Le d\'efaut de la conjecture
d\'epend largement du fait que la fonction de Green actuelle n'est
pas extr\'emale. Il est alors tentant de modifier la conjecture de
mani\`ere \`a la restreindre aux fonctions de Green qui sont
extr\'emales. Pourtant, nous estimons qu'il n'y a pas suffisamment
de support pour faire une conjecture d'aucune fa\c con.

\section{La fonction de Green \`a plusieurs p\^oles}
Commen\c cons avec un lemme connu \cite{johan},\cite{christer}.
\begin{lemma}\label{lokal} Soient $G$ une boule de $\mathbb{C}^n$,
centr\'ee \`a l'origine, et $u$ une fonction plurisousharmonique
n\'egative dans $G$. On pose $$\Psi_u(z,r):=\frac{1}{\log
r}\sup_{\xi\in rD}u(\xi z),$$ pour  $z\in (1/r)G$ et $0<r<1$.
Alors
\newline 1) pour tout $z$ fix\'e, $\Psi_u(z,r)$ est une fonction croissante
de $r$ (sur son domaine de d\'efinition), donc la limite
$$\Psi_u(z)=\lim_{r\downarrow 0}\Psi_u(z,r)$$ existe,\newline 2)
pour $z$ fix\'e, ou bien $u_z:\xi\mapsto u(\xi z)$ est
identiquement $-\infty$, ou bien elle est sousharmonique et on a
$\Psi_u(z)=\triangle u_z(\{0\})$, donc $\Psi_u(z)$ est lin\'eaire
en $u$,\newline 3) il existe une constante $\alpha\geq 0$ (le
nombre de Lelong de $u$ \`a l'origine) et un ensemble pluripolaire
$E\subset\mathbb{C}^n$, tels que $\Psi_u(z)\equiv\alpha$ si
$z\notin E$ et $\Psi_u(z)>\alpha$ si $z\in E$,\newline 4) pour
chaque $z\in G$ et $|\xi|<1$ on a $u(\xi z)\leq\alpha\log|\xi|$.
\end{lemma}
\begin{proof}
Sans perte de g\'en\'eralit\'e, on peut supposer que $G=B$, la boule
unit\'e.\newline 1) et 2) sont
des r\'esultats de th\'eorie du potentiel classique, cf. par exemple
\cite{ransford}, p 46 et
78.\newline 3) Il r\'esulte de 1) que $$-\Psi_u(z)=\lim_{r\downarrow
0}-\Psi_u(z,r)=\sup_{R>r>0}-\Psi_u(z,r)=\sup_{R>r>0}\frac{1}{-\log
r}\sup_{|\xi|=r}u(\xi z),$$ pour
tout $z\in\mathbb{C}^n$ et $R$ assez petit. Par cons\'equent,
$(-\Psi_u)^*\in PSH(\mathbb{C}^n)$ et
comme elle est n\'egative elle est constante ($=-\alpha$). On pose
$E=\{z\in\mathbb{C}^n:
(-\Psi_u)^*(z)\neq -\Psi_u(z)\}$. Alors $E$ et $\alpha$ ont les propri\'et\'es
souhait\'ees.\newline 4) Gr\^ace \`a la monotonie on a $$\frac{1}{\log
r}\sup_{|\xi|=r}u(\xi
z)=\Psi_u(z,r)\geq\Psi_u(z)\geq\alpha,$$ donc $$\sup_{|\xi|=r}u(\xi
z)\leq\alpha\log|\xi|,$$ ce qui
ach\`eve la preuve.
\end{proof}

Soient $\Omega$ un domaine {\em hyperconvexe born\'e} dans $\mathbb{C}^n$,
et $A$ comme dans
l'introduction. Il est commode d'introduire un ordre partiel sur
$(\mathbb{R}^+)^k$ d\'efini par
$\mu=(\mu_1,\ldots,\mu_k)\leq(\nu_1,\ldots,\nu_k)=\nu$ si
$\mu_j\leq\nu_j,j=1,\ldots,k$. Les
p\^oles \'etant fix\'es, on note $g_\nu$ la fonction de Green de poids
$\nu_j$ en $w_j$,
$j=1,\ldots,k$, et $\nu\in(\mathbb{R}^+)^k$ (voir la d\'efinition dans
l'introduction). On a le
th\'eor\`eme suivant.

\begin{theorem}\cite{Lelong}\label{lelong} La fonction $g_\nu$ est l'unique
solution
du probl\`eme de Dirichlet suivant :
\[
\left\{
\begin{array}{l}
u\in C(\overline{\Omega }\setminus A )\cap PSH(\Omega ) \\
(dd^cu)^n=0\text{ dans }\Omega \setminus
A \\ u(z)-\nu_j\log \left| z-w_j\right| =O(1)\text{ si }z\rightarrow
w_j\text{ pour tout }j \\
u(z)\rightarrow 0\text{ si }z\rightarrow \partial \Omega \text{.}
\end{array}
\right.
\]
\end{theorem}

La proposition suivante est une g\'en\'eralisation d'un th\'eor\`eme de
\cite{johan}.

\begin{proposition}\label{green}
Supposons que $g_\nu=\varphi_1+\varphi_2$. Alors, il existe
$\lambda\in(\mathbb{R}^+)^k$ tel que
$\lambda\leq\nu$, $\varphi_1=g_\lambda$ et $\varphi_2=g_{\nu-\lambda}$.
\end{proposition}

\begin{proof}
Notons $\lambda_j$ le nombre de Lelong de $\varphi_1$ et $\mu_j$ celui de
$\varphi_2$ en $w_j$. On
commence par faire une \'etude locale en chaque p\^ole $w_j$. Sans perte de
g\'en\'eralit\'e on
peut supposer que $j=1$ et $w_1=0$. En utilisant les 2) et 3) du lemme
\ref{lokal} on obtient
$\nu_1\equiv\Psi_{g_\nu}(z)=\Psi_{\varphi_1}(z)+\Psi_{\varphi_2}(z)
\geq\lambda_1 + \mu_1$, avec
\'egalit\'e presque partout. Il s'ensuit que $\nu_1=\lambda_1+\mu_1$ et que
les ensembles
exceptionnels de $\varphi_1$ et $\varphi_2$ sont en r\'ealit\'e vides.
D'apr\`es le 4) du lemme
\ref{lokal}, on a $\varphi_1(\xi z)\leq\lambda_1\log|\xi|$ si $|\xi|<1$. Donc
$a(z)=\varphi_1(z)-\lambda_1\log|z|$ est sup\'erieurement born\'ee au
voisinage de $0$, et de la
m\^eme mani\`ere nous trouvons que $b(z)=\varphi_2(z)-\mu_1\log|z|$ y est
sup\'erieurement
born\'ee. Mais leur somme est \'egale \`a $g_\nu(z)-\nu_1\log|z|$ qui est
inf\'erieurement born\'ee
au voisinage de $0$, alors $a$ et $b$ le sont aussi. Nous avons montr\'e
que $\varphi_1$ et
$\varphi_2$ ont le comportement souhait\'e en chaque p\^ole, et que
$\nu=\lambda+\mu$.

On a de plus, $$0=(dd^c g)^n=(dd^c\varphi_1)^n+(dd^c
\varphi_2)^n+\ldots$$ en dehors de $A$, o\`u les termes restants
sont positifs, ce qui montre que $(dd^c\varphi_i)^2=0$, $i=1,2$.
Le th\'eor\`eme \ref{lelong} donne maintenant
$\varphi_1=g_\lambda$, et $\varphi_2=g_\mu=g_{\nu-\lambda}$.
\end{proof}

\begin{corollary}
Pour d\'emontrer que $g_\nu$ est extr\'emale, il suffit de
d\'emontrer que $g_\nu=g_\lambda+g_\mu$ entra\^\i ne que le
vecteur $\lambda$ est proportionnel \`a $\nu$. En particulier,
dans le cas d'un seul p\^ole, $g_\nu$ est toujours extr\'emale.
\end{corollary}
\begin{proof}
Soit donn\'ee une d\'ecomposition $g_\nu=\varphi_1+\varphi_2$. D'apr\`es la
proposition, on a en
r\'ealit\'e $g_\nu=g_\lambda+g_\mu$. Par hypoth\`ese on sait que
$\lambda=c\nu$ o\`u $0\leq c\leq
1$. Comme $g_{c\nu}=cg_\nu$, $g_\nu$ est extr\'emale.
\end{proof}

Supposons que les p\^oles $A=\{w_1\ldots,w_k\}$ et les poids
$\nu=(\nu_1,\ldots,\nu_k)$ sont
fix\'es. Soit $P$ le sous-ensemble de $(\mathbb{R}^+)^k$ d\'efini par
$$\mu=(\mu_1,\ldots,\mu_k)\in
P \Leftrightarrow g_\nu=g_\mu+ g_{\nu-\mu}.$$ Alors, $\{c\nu:\,0\leq c\leq
1\}\subset P$ avec
\'egalit\'e si et seulement si $g_\nu$ est extr\'emal.
\begin{proposition}
L'ensemble $P$ est convexe.
\end{proposition}
\begin{proof}
Supposons que $\mu,\lambda \in P$. Alors,
$$g_\nu=ag_\nu+(1-a)g_\nu=a[g_\mu+g_{\nu-\mu}]+(1-a)
[g_\lambda+g_{\nu-\lambda}]=
[ag_\mu+(1-a)g_\lambda]+.....$$ Maintenant il suffit de montrer que
$ag_\mu+(1-a)g_\lambda=g_{a\mu+(1-a)\lambda}$. On a en dehors de $A$ :
\begin{equation}
\begin{split}
0&=(dd^c g_\nu)^n=\left(dd^c[g_\mu+g_{\nu-\mu}]\right)^k\wedge
\left(dd^c[g_\lambda+g_{\nu-\lambda}]\right)^{n-k}\\ &=(dd^c g_\mu)^k\wedge
(dd^cg_\lambda)^{n-k}+\ldots,
\end{split}
\end{equation}

o\`u les termes restants sont positifs. Il s'ensuit que $$(dd^c
[ag_\mu+(1-a)g_\lambda])^n=\sum_{k=0}^n
\binom{n}{k}a^k(1-a)^{n-k}(dd^cg_\mu)^k\wedge
(dd^cg_\lambda)^{n-k}=0$$ en dehors de $A$. Comme
$ag_\mu+(1-a)g_\lambda$ a le comportement correct en chaque
p\^ole, le th\'eor\`eme \ref{lelong} termine la preuve.
\end{proof}

Passons au cas de la boule unit\'e avec $k=2$ et des poids
\'egaux. Comme nous avons signal\' e dans l'introduction, Coman
\cite{dan} a calcul\'e la fonction de Green correspondante.
Rappellons une partie importante de ce calcul. Apr\`es un
automorphisme appropri\'e, on peut supposer que les p\^oles se
situent sym\'etriquement en $w_1=(-\beta,0)$ et $w_2=(\beta,0)$,
o\`u $\beta\in (0,1)$. Soient $z=(0,\gamma)$, et $S$ l'ensemble
des paires $(s,t)\in (0,1)\times D$ telles qu'il existe un disque
analytique $f:D\to B$ avec $f(0)=z$, $f(s)=w_1$, et $f(t)=w_2$.
Alors, $$S=\{(s,t)\in (0,1)\times D:s\neq t,
s^2>c,|t|^2>c,E(s,t)\geq 0\},$$ o\`u $c$ et $d$ sont des
constantes (qui d\'ependent de $\beta$ et $\gamma$) et
$$E(s,t)=(s^2-c)(|t|^2-c)|1-st|^2-(1-s^2)(1-|t|^2)|st+d|^2.$$
Gr\^ace \`a la sym\'etrie, la fonction de Lempert est r\'ealis\'ee
par un disque correspondant \`a un point de $S$ avec $t$ r\'eel.

En utilisant le calcul de Coman, nous montrerons que la fonction de Green
correspondante est
extr\'emale.
\begin{theorem}
Soient $\Omega=B$, la boule unit\'e dans $\mathbb{C}^n$, $n\geq 2$, et
$w_1$ et $w_2$ deux p\^oles
donn\'es, avec $\nu_1=\nu_2=1$. Alors, $g_{(1,1)}$ est extr\'emale.
\end{theorem}
\begin{proof}
Commen\c cons avec le cas $n=2$. Sans perte de g\'en\'eralit\'e on peut
supposer que les p\^oles se
situent sym\'etriquement. Nous sommes alors dans le cas d\'ecrit ci-dessus.
Supposons que
$g_{(1,1)}=g_{(p,q)}+g_{(1-p,1-q)}$. Nous avons alors
\begin{equation}
\begin{split}\label{dan1}
g_{(1,1)}=& g_{(p,q)}+g_{(1-p,1-q)} \leq \delta_{(p,q)}+\delta_{(1-p,1-q)}
\\=&
\inf\{p\log|\zeta_1|+q\log|\zeta_2|\}+\inf\{(1-p)\log|\zeta_1|+(1-q)\log|
\zeta_2 |\}
\\ \leq &
\inf\{\log|\zeta_1|+\log|\zeta_2|\}=\delta_{(1,1)}=g_{(1,1)},
\end{split}
\end{equation}
o\`u la derni\`ere \'egalit\'e est le th\'eor\`eme de Coman et toutes les
bornes inf\'e\-rieures sont
prises sur la m\^eme famille de disques analytiques. En particulier, le
disque extr\'emal de
$\delta_{(1,1)}$ est extr\'emal pour $\delta_{(p,q)}$ et
$\delta_{(1-p,1-q)}$ aussi. Donc, les
fonctions $S\to\mathbb{R}:$ $$(s,t)\mapsto s|t|,(s,t)\mapsto s^p
|t|^q,(s,t)\mapsto s^{1-p}
|t|^{1-q}$$ sont minimales au m\^eme point $a\in\partial S$. Coman a
montr\'e que $E(a)=0$ tandis
que $\partial S\in C^1$ au voisinage de $a$. Alors les gradients de ces
 trois fonctions sont proportionnels, et $p=q$.
La preuve est finie dans le cas $n=2$.

Passons au cas g\'en\'eral. On suppose que les p\^oles se situent
sym\'etri\-quement dans le disque $\{z\in B_n : z_2=\ldots
=z_n=0\}$, et on note $z=(z_1,z_2,\ldots,z_n)=(z_1,z')$. Alors, la
fonction de Green $g^n_{(1,1)}(z)$ n'est que la fonction de Green
de la boule de dimension 2, $g^2_{(1,1)}$, \'evalu\'ee au point
$(z_1,||z'||)$, cf. \cite{dan}, corollaire 4.6.3. On note
$\delta^n$ et $\delta^2$ les fonctions de Lempert correspondantes.
Soit $0\leq z_2<1$ et $U$ une rotation unitaire dans
$\mathbb{C}^{n-1}$ envoyant $(z_2,0,\ldots,0)$ sur $z'$. Si
$\varphi=(\varphi_1,\varphi_2)$ est un disque analytique
appartenant \`a la famille qui d\'efinit
$\delta^2_{(p,q)}(z_1,z_2)$, on peut fabriquer un disque
$\tilde{\varphi}=(\tilde{\varphi}_1,\tilde{\varphi}_2)$ qui
appartient \`a la famille qui d\'efinit $\delta^n_{(p,q)}(z_1,z')$
en posant $\tilde{\varphi}_1=\varphi_1$ et
$\tilde{\varphi}_2=U\circ\varphi_2$. \'Evidemment, on peut aussi
faire l'inverse. Ceci montre que
$\delta^n_{(p,q)}(z_1,z')=\delta^2_{(p,q)}(z_1,||z'||)$.
Maintenant on peut faire un calcul analogue au pr\'ec\'edent
(\ref{dan1}) :
\begin{equation}
\begin{split}
\delta^n_{(1,1)}(z_1,z')=&\delta^2_{(1,1)}(z_1,||z'||)=g^2_{(1,1)}(z_1,||z'|
|) =g^n_{(1,1)}(z_1,z')
\\=& g^n_{(p,q)}(z_1,z')+g^n_{(1-p,1-q)}(z_1,z')
\\ \leq &
\delta^n_{(p,q)}(z_1,z')+\delta^n_{(1-p,1-q)}(z_1,z')\leq\delta^n_{(1,1)}(z_
1,z' )
\end{split}
\end{equation}
On trouve une contradiction comme plus haut.
\end{proof}

Soient maintenant $D\times D$ le bidisque de $\mathbb{C}^2$, et $a_i\in D$,
$i=1\ldots,k$. On note
$T_i(z_1)=\log|\frac{z_1-a_i}{1-\bar{a}_i z_1}|$ la fonction de Green du
disque unit\'e. Alors la
fonction de Green de p\^oles $w_i=(a_i,0)$ et poids
$\mathbf{1}=(1,\ldots,1)$ est connue
\cite{magnus} :
\begin{equation}
g_\mathbf{1}(z)=\max\{\sum_{i=1}^k T_i(z_1), \log|z_2|\}.\label{green1}
\end{equation}
G\'en\'eralisons cette formule. On garde la position des p\^oles, mais on
autorise des poids
diff\'erents. Sans perte de g\'en\'eralit\'e on peut supposer que
$1=\nu_1\geq\nu_2\geq\ldots\geq\nu_k$.
\begin{theorem} \label{bidisk} Dans la situation d\'ecrite, si au moins
deux poids sont
diff\'erents, la fonction de Green est la suivante :
\begin{equation}
g_\nu=
g_{(1,\nu_2,\ldots,\nu_k)}=
\nu_k h_k(z)+\sum_{j=1}^{k-1} (\nu_j-\nu_{j+1})h_j(z),\label{green2}
\end{equation}
o\`u
$h_j=\max\{T_1(z_1)+\ldots+T_j(z_1),\log|z_2|\}$ est la fonction de Green
de $D\times D$ avec le
poids $1$ en $w_1,\ldots,w_j$.

Par cons\'equent, elle n'est pas extr\'emale. En revanche, si tous
les poids sont \'egaux (cf. la formule (\ref{green1})) elle est
extr\'emale.
\end{theorem}
\begin{proof}
Soit $b(z)$ la somme de (\ref{green2}). Commen\c cons par montrer que
$g=b$. Il suffit de
v\'erifier les quatre conditions du th\'eor\`eme \ref{lelong}. Les
premi\`ere et quatri\`eme
conditions sont triviales. On v\'erifie la troisi\`eme condition. Au
voisinage de $w_l$, la
fonction $h_j$ peut s'\'ecrire $h_j(z)=\log||z-w_l||+O(1)$ si $j\geq l$, et
elle y est born\'ee si
$j<l$. Il s'ensuit que
$$b(z)=\left(\nu_k+\sum_{j=l}^{k-1}(\nu_j-\nu_{j+1})\right)\log||z-w_l||+O(1)=
\nu_l\log||z-w_l||+O(1)$$ au voisinage de $w_l$. Il nous reste \`a
v\'erifier que $b$ est maximale
en dehors des p\^oles. On constate que, pour tout point except\'e les
p\^oles, il existe un
voisinage $U$ o\`u toutes les $h_j$, sauf peut-\^etre une, sont
pluriharmoniques. Donc, dans $U$,
la masse de Monge-Amp\`ere de $b$ est donn\'ee par la fonction
exceptionnelle et comme toutes les
$h_j$ sont maximales on conclut que $g=b$.

La formule (\ref{green2}) implique \'evidemment la
non-extr\'emalit\'e de $g_\nu$ dans le cas o\`u les poids sont
distincts. Finalement, supposons que
$g_\mathbf{1}=g_\nu+g_{\mathbf{1}-\nu}$, o\`u les $\nu_j$ sont
dans un ordre d\'ecroissant. On peut \'ecrire en utilisant
(\ref{green2}) $$g_\nu=\nu_k
h_k+\sum_{j=1}^{k-1}(\nu_j-\nu_{j+1})h_j.$$ De m\^eme on a
$$g_{\mathbf{1}-\nu}=(1-\nu_1)h_k+\sum_{j=1}^{k-1}(\nu_j-\nu_{j+1})h'_j,$$
o\`u $h'_j$ est la fonction de Green associ\'ee \`a l'ensemble des
p\^oles,
 compl\'ementaire
de celui de $h_j$. Cela entra\^\i ne que la fonction $h_j+h'_j$ a
des p\^oles de poids $ \mathbf{1}$ en tout $w_j$, $1\leq j\leq k$.
On trouve par substitution
$$(\nu_1-\nu_k)h_k=\sum_{j=1}^{k-1}(\nu_j-\nu_{j+1})(h_j +h'_j).$$
Comme pour $j<k$ on a $h_j+h'_j\le h_k$ avec in\'egalit\'e stricte
quelque part dans $D\times D$, on conclut que tous les $\nu_k$
sont \'egales, et le th\'eor\`eme est d\'emontr\'e.
\end{proof}

\remark Nous donnons une autre forme de la fonction de Green \`a
poids diff\'erents\ :
$$g_\nu=\max\{u_1(z),\ldots,u_k(z),\log|z_2|\} $$ o\`u
$$u_1(z)=\sum_{i=1}^k\nu_i T_i(z_1),$$
$$u_j(z)=\nu_j\log|z_2|+\sum_{i=1}^{j-1}(\nu_i-\nu_j)T_i(z_1),\quad
2\leq j\leq k.$$ Nous esquissons la preuve. Soit $a(z)$ le max. On
a $$U_1:=\{z:a(z)=u_1(z)\}=\{z:\log|z_2|\leq\sum_{i=1}^k
T_i(z_1)\},$$ $$U_j:=\{z:a(z)=u_j(z)\}=\{z:\sum_{i=1}^j
T_i(z_1)\leq\log|z_2|\leq\sum_{i=1}^{j-1} T_i(z_1)\},$$ pour
$j=2,\ldots,k$, et $$V:=\{z:a(z)=\log|z_2|\}=\{z:\log|z_2|\geq
T_1(z_1)\}.$$ Remarquons que ces ensembles sont invariants par
changements des poids, tant que {\it l'ordre} des poids est
inchang\'e.

Dans $U_1$, on a $h_j(z)=T_1(z_1)+\ldots+T_j(z_1)$ et dans $V$ on a
$h_j(z)=\log|z_2|$ pour tout
$j$. Dans $U_l$, $l=2,\ldots,k$, on a $h_j(z)=T_1(z_1)+\ldots+T_j(z_1)$ si
$1\leq j\leq l-1$ et
$h_j(z)=\log|z_2|$ si $j\geq l$. En utilisant cela et en consid\'erant
chaque ensemble
s\'epar\'ement, on peut v\'erifier que $a=b$.

\remark Le th\'eor\`eme reste vrai pour le polydisque dans
$\mathbb{C}^n$, si tous les p\^oles se situent dans le disque
$\{z_2=\ldots=z_n=0\}$. Dans ce cas, il faut remplacer $\log|z_2|$
par $v(z)=\max\{\log|z_2|,\ldots,\log|z_n|\}$.

\section{La fonction $\max\{\log|z|,-1\}$ dans la boule unit\'e}
Soit $u(z)=\max\{\log|z|,-1\}$ la fonction relative extr\'emale du compact
$B_1=\{z:|z|\leq 1/e\}$
dans la boule unit\'e $B$ de $\mathbb{C}^2$. On note $B_2=B\setminus B_1$.
On consid\`ere la
question de l'extr\'emalit\'e de $u$.
\begin{theorem} La fonction $u$ est extr\'emale.
\end{theorem}
\begin{proof}
Supposons que $u=\varphi_1+\varphi_2$, o\`u $\varphi_i$ est une fonction
plurisous-harmonique
n\'egative dans la boule $B$, $i=1,2$. Remarquons que les fonctions
$\varphi_i$ sont sup\'erieurement
continues et leur somme est continue : elles sont alors continues.
De plus, $-1<\varphi_i<0$.

Dans $B_1$, la fonction $u$ est pluriharmonique. Dans $B_2$, elle
est harmonique sur chaque droite passant par l'origine. Plus
pr\'ecisement, pour chaque $q\in\mathbb{C}$ fix\'e, $z_1\mapsto
u(z_1,qz_1)$ est harmonique dans la couronne $1/(e
\sqrt{1+|q|^2})<|z_1|<1/\sqrt{1+|q|^2}$. (De plus, $z_2\mapsto
u(0,z_2)$ est harmonique dans la couronne $1/e<|z_2|<1$, ce qui
correspond \`a $q=\infty$.) Par cons\'equent, $\varphi_1$ et
$\varphi_2$ ont \'egalement toutes les propri\'et\'es
mentionn\'ees.

On va d\'emontrer le th\'eor\`eme sous une condition
suppl\'ementaire, que $\varphi_1$ (et donc $\varphi_2$) ne
d\'epend que de $z_1$ et $|z_2|$, qui sera enlev\'ee apr\`es.
D'abord $\varphi_1$ ne d\'epend que de $z_1$ dans $B_1$. Pour voir
cela, on fixe $z_1$. Alors $\varphi_1$ est une fonction harmonique
de $z_2$ dans un disque centr\'e en $0$, et ne d\'epend que de
$|z_2|$, donc elle est constante. Il s'ensuit que la fonction
suivante est bien d\'efinie : $v(z_1)=\varphi_1|_{B_1}(z_1,z_2)$
pour $|z_1|<1/e$. On pose aussi $V(z)=\varphi_1(z)$ dans $B_2$. La
fonction $v$ est harmonique et continue jusqu'au bord. La fonction
$V$ est confondue avec $v$ sur $|z|=1/e$, elle est continue
jusqu'au bord, elle s'annule sur $\partial B$. De plus elle est
harmonique sur chaque droite complexe passant par l'origine. Ceci
montre que pour chaque $v$ donn\'ee, $V$ est unique si elle
existe; s'il y avait deux telles fonctions, on consid\'ererait
leur diff\'erence sur les droites.

La fonction $v$ admet la repr\'esentation
$$v(z_1)=\sum_{n=-\infty}^\infty c_n (er)^{|n|} e^{in\theta}$$
dans le disque $|z_1|\leq 1/e$. Ici, $z_1=re^{i\theta}$, et
$$c_n=\frac{1}{2\pi}\int_0^{2\pi}v(e^{-1+it})e^{-int}dt.$$
Note que $c_0=v(0)=\varphi_1(0)\in (-1,0)$. On
cherche ensuite une fonction $H_t(w)$, $0\leq t<1/e$, qui soit
harmonique dans la couronne $t<|w|<et$, avec les valeurs au bord
$H_t(w)=0$ sur le cercle $|w|=et$ et $H_t(w)=v(w)$ sur $|w|=t$. On
v\'erifie sans peine que
$$H_t(w)=-c_0\log|\frac{w}{et}|+\sum_1^\infty\left(c_n
(ew)^n+c_{-n}(e\bar{w})^n\right)\frac{(et/|w|)^{2n}-1}{e^{2n}-1}$$
est la solution unique de ce probl\`eme.

Consid\'erons maintenant, pour $z_2/z_1$ fix\'e, la fonction $w\mapsto
V(w,(z_2/z_1)w)$. Si l'on
pose $t=|z_1|/(e|z|)$, elle se confond avec $H_t$. Par cons\'equent,
\begin{equation}
\begin{split}
V(w,(z_2/z_1)w)=&-c_0\log|\frac{w|z|}{z_1}|\\+&\sum_1^\infty\left(c_n
(ew)^n+c_{-n}(e\bar{w})^n\right)\frac{(|z_1|/(|z||w|))^{2n}-1}{e^{2n}-1}
.
\end{split}
\end{equation}

En particulier, en posant $w=z_1$ on trouve
$$V(z_1,z_2)=-c_0\log|z|+\sum_1^\infty\left(c_n
(ez_1)^n+c_{-n}(e\bar{z}_1)^n\right)\frac{(1/|z|)^{2n}-1}{e^{2n}-1}.$$
Pour abr\'eger, rempla\c cons
$c_n e^{|n|}/(e^{2|n|}-1)$ par $b_n$, $n\neq 0$, et $-c_0$ par
$a$; la formule se ram\`ene \`a
        $$V(z_1,z_2)=a\log|z|+\sum_1^\infty\left(b_n
        z_1^n+b_{-n}\bar{z}_1^n\right)\left((1/|z|)^{2n}-1\right).$$
Comme $V$ a des valeurs r\'eelles, on
        obtient $b_{-n}=\bar{b}_n$, de sorte que $V(z)=a\log|z|+\real
        f(z_1/|z|^2)-\real f(z_1)$ o\`u $f$
        est une fonction holomorphe dans le disque de rayon $e$ (car
$1/e<|z|<1$, on a
        $0\leq|z_1|/|z|^2<e$). Posons $g(z)=f(1/z)$ : $g$ est holomorphe au
dehors du disque de rayon
        $1/e$, et $V(z)=a\log|z|+\real g(\bar{z}_1+z_2\bar{z}_2/z_1)-\real
g(1/z_1)$.

  On calcule ensuite le signe du d\'eterminant $(dd^ch)^2$ de la matrice
        Hessienne de
        $h(z):=V(z)-a\log|z|=\real g(\bar{z}_1+z_2\bar{z}_2/z_1)-\real
g(1/z_1)$.
        Le dernier terme est
        pluriharmonique. Il suffit alors de consid\'erer
        $g(\bar{z}_1+z_2\bar{z}_2/z_1)$. Il est commode de
        faire le changement de coordonn\'ees suivant: $w_1=z_2/z_1,w_2=z_1$.
Comme il est holomorphe, le
        signe du d\'eterminant ne changera pas, et
        \begin{equation}
        \begin{split}
        2\real g(\bar{z}_1+z_2\bar{z}_2/z_1)=&2\real
g(\bar{w}_2(1+w_1\bar{w}_1))\\=&
        g(\bar{w}_2(1+w_1\bar{w}_1))+g(w_2(1+w_1\bar{w}_1))=:\tilde{h}(w).
        \end{split}
        \end{equation}
        La derni\`ere expression \'etant harmonique en $w_2$, il est
superflu de calculer
        $\partial^2\tilde{h}/\partial w_1\partial\bar{w}_1$, et le
d\'eterminant vaut
        $$-\left|\frac{\partial^2 \tilde{h}}{\partial
        w_1\partial\bar{w}_2}\right|^2=-\left|\bar{w}_1(D_1+D_2\bar{w}_2
        (1+w_1\bar{w}_1))\right|^2,$$
o\`u $D_j$ note la d\'eriv\'ee $j$-i\`eme de $g$ \'evalu\'ee au point
        $\bar{w}_2(1+w_1\bar{w}_1)$. On
        conclut que, ou bien le d\'eterminant de la matrice Hessienne de $h$
est n\'egatif quelque part, ou
        bien $D_1+D_2\bar{w}_2 (1+w_1\bar{w}_1)$ s'annule partout.

        Si le d\'eterminant est n\'egatif en un point, nous allons
        d\'eduire une contradiction. En effet, dans $B_2$,
        $$(dd^c\varphi_1)^2=(dd^c V(z))^2=(dd^c h(z))^2+2a\,
        dd^c\log|z|\wedge dd^c h(z),$$ et
        $$(dd^c\varphi_2)^2=(dd^c(\log|z|-V(z)))^2=(dd^c
        h(z))^2-2(1-a)dd^c\log|z|\wedge dd^c h(z).$$
Car $a=-c_0\in(0,1)$, il r\'esulte que, au
        point o\`u $(dd^ch)^2$ est n\'egatif, ou bien
        $(dd^c\varphi_1)^2<0$, ou bien $(dd^c\varphi_2)^2<0$, ce qui est
        contradictoire.

        Si, d'autre part, $D_1+D_2\bar{w}_2 (1+w_1\bar{w}_1)$ s'annule
partout, on
        voit, en rempla\c cant
        $\bar{w}_2 (1+w_1\bar{w}_1)$ par $z$, que $g'(z)+zg''(z)=0$ partout.
Mais
        cela est impossible
        lorsque $g$ est holomorphe au voisinage de l'infini, sauf si $g$ est
        constante. Alors $f\equiv 0$,
        $V(z)=a\log|z|$ et par continuit\'e $\varphi_1(z)=au(z)$. Ceci
ach\`eve la
        d\'emonstration sous la
        condition suppl\'ementaire.

        Passons au cas g\'en\'eral. On d\'efinit
        \begin{equation}
        \Phi_i=\frac{1}{2\pi}\int_0^{2\pi}\varphi_i(z_1,z_2
e^{i\theta})d\theta,\quad
        i=1,2.\label{integrale}
        \end{equation}
        Alors $g=\Phi_1+\Phi_2$, et $\Phi_i$ ne d\'epend que de $z_1$ et
$|z_2|$,
        $i=1,2$, donc le cas
        sp\'ecial montre qu'il existe une constante $a\geq 0$ telle que
        $\Phi_1(z)=au(z)$ pour tout $z\in
        B$. Or sur la droite complexe $z_2=0$, on a
$\Phi_1(z)=\varphi_1(z)$, ce
        qui implique que
        $\varphi_1(z)=au(z)$ sur cette droite.

        Finalement, si $d$ est une droite complexe quelconque passant par
        l'origine, il existe une
        transformation unitaire $R_d$, telle que $R_d^{-1}(d)=\{z_2=0\}$.
Comme $u$
        est invariante par
        cette transformation, on peut \'evidemment remplacer $\varphi_i$ par
        $\varphi_i\circ R_d$ dans
        l'argument pr\'ec\'edent. Ceci montre que $\varphi_1(z)=a_d u(z)$
sur $d$, o\`u $a_d\geq 0$.
        L'origine \'etant un point commun \`a toutes les droites, on obtient
le r\'esultat escompt\'e.
        \end{proof}

        \remark Le th\'eor\`eme reste vrai pour $n\geq 3$, avec presque la
        m\^eme d\'emon\-stration. Il faut d'abord supposer que les
        composantes ne d\'ependent que de $z_1$ et $||z'||$, o\`u
        $z'=(z_2,\ldots,z_n)$. On conclut qu'elles sont proportionnelles
        \`a $u$. Ensuite, on remplace la formule (\ref{integrale}) par
        
$$\Phi_i=\frac{1}{(2\pi)^{n-1}}\int_{[0,2\pi]^{n-1}}\varphi_i(z_1,z_2
        e^{i\theta_2},\ldots,z_n e^{i\theta_n})d\theta_2\cdots
        d\theta_n,\quad i=1,2,$$ etc.

\section{Fonctions extr\'emales qui ne s'annulent pas sur tout le bord}
        Si $\sup_{z\in\Omega}u(z)=c<0$, alors $u$ n'est pas extr\'emale
        puisqu'elle peut \^etre d\'ecompos\'ee : $u(z)=h/2+(u(z)-h/2)$,
o\`u $h$ est une fonction plurisousharmonique avec $c<h<0$,
 qui n'est pas une multiple de $u$. Il
        s'ensuit que toute fonction extr\'emale, continue jusqu'au bord,
        s'annule quelque part au bord.

        D'autre part,  il est facile de donner un exemple d'une fonction
        extr\'emale qui est n\'egative sur une grande partie du bord. On
        prend simplement la fonction $\log|z_1|$ dans le bidisque.
        Supposons que $\log|z_1|=\varphi_1(z)+\varphi_2(z)$. Si on fixe
        $z_2$, on a
        $\log|\cdot|=\varphi_1(\cdot,z_2)+\varphi_2(\cdot,z_2)$. Comme le
        logarithme est extr\'emal dans le disque unit\'e, on conclut que
        $\varphi_1(z)=c(z_2)\log|z_1|$. Ensuite, en fixant $z_1$, on
        trouve que $c(z_2)$ est harmonique. On calcule $(dd^c
        u)^2(z)=-|\partial\log|z_1|/\partial z_1|^2 |\partial
        c(z_2)/\partial z_2|^2\leq 0$. Comme $u$ est plurisousharmonique
        l'expression s'annule partout. En cons\'equence, $c$ est
        constante.

        Si $\Omega$ est un domaine B-r\'egulier (par exemple la boule), le
        probl\`eme de Dirichlet de l'\'equation de Monge-Amp\`ere a une
        solution pour toute fonction $f$ continue sur le bord. On peut se
        demander pour quelles fonctions $f$ la solution est extr\'emale.
        Il faut que l'ensemble o\`u $f$ s'annule soit suffisamment large,
        mais le probl\`eme reste myst\'erieux.

        \section{La conjecture de Coman}

        Les d\'efinitions utilis\'ees dans cette section se trouvent dans
        l'introduction. La fonction de Green dans le bidisque, que nous
        avons calcul\'ee plus haut, donne un contre-exemple \`a la
        conjecture de Coman. Pour montrer cela nous utiliserons les
        th\'eor\`emes \ref{Frank1} et \ref{bidisk} et le r\'esultat
        suivant.
\begin{theorem}\cite{Frank}{\em ,  th\'eor\`eme 2.4.}\label{Frank2}
        Soient $\Omega$ un domaine born\'e convexe de $\mathbb{C}^n$ et $A$
comme plus haut. Alors $\delta^A(z)=\delta(z,A)$ pour tout $z\in\Omega$.
        \end{theorem}

        Dans le reste de l'article nous fixons deux p\^oles $(a,0)$ et
        $(b,0)$ dans le bidisque, $a\neq b$ et $a,b\neq 0$, et un point
        $z=(0,\gamma)$, tel que $|ab|<|\gamma|<\min\{|a|,|b|\}$. On note
        $g_{p,q}$ la fonction de Green function avec poids $p$ en $(a,0)$
        et $q$ en $(b,0)$, \'evalu\'ee en $z$, et de m\^eme
        $\delta_{p,q}$. Remarquons que, d'apr\`es le th\'eor\`eme
        \ref{bidisk}, $g_{1,1}=\log|\gamma|$, $g_{1,0}=\log|a|$, et
        $g_{2,1}=g_{1,1}+g_{1,0}=\log|\gamma|+\log|a|$.

        \smallskip
        Maintenant nous pouvons d\'emontrer notre th\'eor\`eme.

        {\it D\'emonstration du th\'eor\`eme \ref{counterexample}.} Il
        suffit de montrer que $\delta_{2,1}>g_{2,1}$. On sait d\`eja que
        $\delta_{2,1}\geq g_{2,1}$. Supposons que l'in\'egalit\'e est une
        \'egalit\'e. Alors
        \begin{equation}
        \begin{split}
        g_{2,1}=& g_{1,1}+g_{1,0} \leq \delta_{1,1}+\delta_{1,0}\\ =&
        \inf\{\log|\zeta_1|+\log|\zeta_2|\}+\inf\{\log|\zeta_1|\}
        \\ \leq &\inf\{2\log|\zeta_1|+\log|\zeta_2|\}=\delta_{2,1}=g_{2,1},
        \end{split}
        \end{equation}
        o\`u toutes les bornes inf\'erieures sont prises sur la m\^eme
        famille $F_z$. Donc, toutes les in\'egalit\'es sont en r\'ealit\'e
        des \'egalit\'es. En utilisant les th\'eor\`emes \ref{Frank1} et
        \ref{Frank2}, la  {\it derni\`ere} borne inf\'erieure est atteinte
        par un disque extr\'emal $f$ qui passe par $(a,0)$ ou $(b,0)$ ou
        tous les deux. Il s'ensuite que $f$ est \'egalement extr\'emal
        pour $\delta_{1,1}$ et $\delta_{1,0}$. Pourtant, cela est
        impossible, d'apr\`es le lemme suivant. La contradiction donne le
        th\'eor\`eme.\endproof

        \begin{lemma}
Si $|ab|<|\gamma|<\min\{|a|,|b|\}$, il n'y a aucun disque extr\'emal commun
        \`a $\delta_{1,1}$ et $\delta_{1,0}$.
        \end{lemma}
        \begin{proof}
        Commen\c cons par caract\'eriser tous les disques extr\'emaux pour
        $\delta_{1,0}$. Soit $f=(f_1,f_2):D\to D\times D$ un tel disque.
        Par d\'efinition, il existe $\zeta_1\in D$ tel que
        $f_1(\zeta_1)=a$ et $f_1(0)=0$. On a $g_{1,0}=\log|a|$, et cela
        est \'egal \`a $\delta_{1,0}$, car cette valeur est atteinte par
        le disque
        
$$\zeta\mapsto(\zeta,\frac{\gamma}{a}\frac{\zeta-a}{1-\bar{a}\zeta}).$$
        Donc $|\zeta_1|=|a|$. En utilisant le lemme de Schwarz on conclut
        qu'un disque passant par $(a,0)$ et $z$ est extr\'emal pour
        $\delta_{1,0}$ si et seulement si il est une rotation dans la
        premi\`ere variable.

        Fixons maintenant un tel disque
        $\zeta\mapsto(\alpha\zeta,f_2(\zeta))$, o\`u $|\alpha|=1$, et
        supposons qu'il est extr\'emal pour $\delta_{1,1}$. D'apr\`es le
        th\'eor\`eme \ref{Frank1} il y a deux possibilit\'es. Soit le
        disque passe par un seul p\^ole, ce qui est fatalement $(a,0)$,
        soit par tous les deux p\^oles. Dans le premier cas, $(a,0)$
        serait l'image de $\zeta_1=a/\alpha$ et on calculerait
        $\delta_{1,1}=\log|a/\alpha|=\log|a|$. D'autre part, si le disque
        passait par les deux p\^oles, ceux-ci seraient les images de
        $\zeta_1=a/\alpha$ et $\zeta_2=b/\alpha$ respectivement. Alors
        $\delta_{1,1}=\log|a/\alpha|+\log|b/\alpha|=\log|ab|$ dans ce cas.

        Pour conclure la preuve, on va montrer que
        $\delta_{1,1}=g_{1,1}=\log|\gamma|$, ce qui exclut les deux
        possibilit\'es. On a $|\gamma|<|a|<|a/b|$, donc $|\gamma b/a|<1$.
        Pareillement on obtient $|\gamma a/b|<1$. Par cons\'equent, on
        peut choisir $\zeta_1\in D$ tel que $\zeta_1^2=\gamma a/b$. Posons
        $\zeta_2=\gamma/\zeta_1$, et $\beta=a/\zeta_1$. Alors $\zeta_2\in
        D$, car $|\zeta_2|^2=|\gamma b/a|$, et $\beta\in D$ car
        $|\beta|^2=|ab/\gamma|$. Maintenant on d\'efinit un disque
        analytique par
        
$$f:\zeta\mapsto\left(\beta\zeta,\frac{\zeta-\zeta_1}{1-\bar{\zeta}_1\zeta}
        \frac{\ \zeta-\zeta_2}{1-\bar{\zeta}_2\zeta}\right).$$ Il est
        facile de v\'erifier qu'il envoie $\zeta_1$ sur $(a,0)$, $\zeta_2$
        sur $(b,0)$, et $0$ sur $(0,\gamma)$. On calcule
        $d(f)=\log|\zeta_1|+\log|\zeta_2|=\log|\gamma|$, ce qui termine la
        preuve du lemme.
        \end{proof}

        \begin{theorem} Soit $A$ comme plus haut. Il existe des domaines
strictement convexes, lisses,
        contenus dans le bidisque, tels que $g(z,A)\not\equiv\delta(z,A)$.
        \end{theorem}
        \begin{proof} Soit $(\Omega_j\subset D\times D)_j$ une suite
croissante
         de domaines convexes, lisses, avec $\cup_j\Omega_j=D\times D$.
        On a $\delta_{\Omega_j}(z,A)\ge \delta(z, A)$, parce que pour
$\Omega_j$
         la borne inf\'erieure est prise par
        rapport \`a une plus petite famille $F_z$.

        Soit $\epsilon>0$ assez petit. Quand $\Omega_j$ est si grand que
        $G_\epsilon=\{z;\, g(z,A)<-\epsilon\}\subset \Omega_j$, alors
        $g(z,A)+\epsilon$ est la fonction de Green de $G_\epsilon$. Par
        cons\'equent $g_{\Omega_j}(z,A)\le g(z,A)+\epsilon$. Ceci montre
        en notre cas le fait connu que, pour $z$ fix\'e, la valeur
        $g_\Omega(z,A)$ varie contin\^ument avec $\Omega$ croissante. On
        a, avec $z=(0,\gamma)$ : $$\delta_{\Omega_j}(z,A)\ge \delta(z,
        A)>g(z,A)\geq g_{\Omega_j}(z,A)-\epsilon.$$ Le th\'eor\`eme
        s'ensuit.
        \end{proof}

        \smallskip
        {\it Remerciements.} Ce travail a commenc\'e lors de notre visite au
        Laboratoire E. Picard,
        Universit\'e Paul Sabatier, Toulouse. Nous remercions ses membres
pour leur
        hospitalit\'e et des
        discussions int\'eressantes. De plus, nous remercions Urban Cegrell,
qui
        nous a inspir\'e pour
        travailler sur ce genre de probl\`emes. Magnus Carlehed a re\c cu
l'aide de
        l'Institut Su\'edois,
        la Fondation de Wenner-Gren et la Fondation en m\'emoire de Lars
Hierta.

        \end{document}